# Monotonically Decreasing the Number of Directed 3-Cycles via Edge-Flips?

David Bom, Florian Unger, and Birgit Vogtenhuber

Technische Universität Graz, Austria
{bom@student.tugraz.at, florian.unger@tugraz.at,
birgit.vogtenhuber@tugraz.at}

**Abstract.** We investigate a combinatorial reconfiguration problem on oriented graphs, where a reconfiguration step (edge-flip) is the inversion of the orientation of a single edge. A recently published conjecture that is relevant to the correctness of a Markov Chain Monte Carlo sampler for directed flag complexes states that any simple oriented graph admits a flip sequence that monotonically decreases the number of directed 3-cycles to zero, and is known to be true for complete oriented graphs.

We show that, in general, this conjecture does not hold. As main tool for disproving the conjecture, we introduce the concept of FBD-graphs (fully blocked digraphs). An FBD-graph is a directed graph that does not contain any directed 1-, 2-, or 3-cycles, and for which any edge-flip creates a directed 3-cycle. We prove that the non-existence of FBD-graphs is a necessary condition for the conjecture to hold and succeed in constructing FBD-graphs. On the other hand, we show that complete graphs, as well as graphs in which every edge is incident to at most two triangles, cannot be fully blocked. In addition to being relevant for determining in which cases the above mentioned sampling process is correct, the concept of FBD-graphs might also be useful for other problems and yields interesting questions for further study.

## 1 Introduction

Combinatorial reconfiguration problems are an ample source of our favourite kind of problem: easy, sometimes even elementary to describe, but often hard to solve. Key characteristics of a combinatorial reconfiguration problem are a (finite) collection of combinatorial configurations and a (small) set of reconfiguration moves, which transform one configuration into another. The moves induce a reconfiguration graph on the configurations, where two configurations are connected via an edge if one can be transformed into the other by a single move. The most fundamental question about such a graph is whether it is connected. In case of a positive answer, typical follow-up questions ask for (bounds on) its diameter or radius, the computational complexity of finding shortest paths, its degree of connectivity, expansion properties, Hamiltonicity, and many more. See the related work section below for some examples. Besides being interesting structures that are worthy to be studied for the joy of problem solving alone, combinatorial reconfiguration problems are often motivated by and/or relevant for applications such as motion planning, sorting, or random sampling, the latter being the case here:

The problem we investigate stems from Markov Chain Monte Carlo (MCMC) sampling of directed flag complexes [28]. MCMC-sampling works by iteratively and stochastically applying many small modifications to the starting point until the point one ends up in is (uniformly) random. In our terms, this is described by a (long) random walk on the reconfiguration graph, whose connectedness is thus an elementary prerequisite to ensure that one samples not only from a subset of valid points.

Thus we consider connectedness of flip graphs on orientations of simple graphs (for short, simple oriented graphs), where the configuration space is the set of all orientations of a given undirected simple graph, and a flip consists of reversing the direction of a single edge. Without any further restriction, the reconfiguration question of connecting any two oriented graphs which share the same underlying simple graph is straightforward - one can simply revert the orientations of all differently oriented edges in an arbitrary order. However, the problem under consideration does have additional restrictions, as motivated in [28], leading them to state the following question (Conjecture 4 in [28]):

*Can any simple, oriented graph be transformed into a directed acyclic graph by a sequence of edge-flips, such that the number of directed 3-cycles monotonically declines?*

They justify stating this conjecture as they proved it for complete oriented graphs, and furthermore verified it on billions of randomly sampled graphs.

**Contribution and Outline.** Our main result lies in refuting this conjecture, which we achieve in two steps: We show that the conjecture can't be true if *FBD-graphs*, which are a novel class of graphs first described here, exist (Section 3); and we construct an FBD-graph (Section 4). We further show that FBD-graphs can't exist in some cases (Section 3), and raise a couple of elementary, but independently interesting open questions (Section 5). In the following we sketch the two steps to disprove the main conjecture:

*Reduction:* What are prohibited edge flips? This question immediately leads to the concept of *sufficiently blocked edges*: these are edges whose inversion would create more 3-cycles than they destroy. A graph consisting solely of unflippable edges thus potentially is a nice building block to refute the conjecture. This motivates the definition of fully blocked digraphs (*FBD-graphs*) which for compliance with the conjecture must be finite, simple, oriented and, in order for a single block per edge to suffice, also 3-cycle free. With that we describe a graph where flipping any edge is illegal, but misses one 3-cycle to form a counterexample by itself.

A *doubly* blocked edge may then be created by gluing two FBD-graphs along an edge, and three doubly blocked edges arranged to form a 3-cycle constitute a counterexample to the conjecture: There is one 3-cycle in the graph, but trying to remove it by flipping *any* of its edges immediately creates two new 3-cycles, increasing the number of 3-cycles overall.

*Construction:* It remains to show that such FBD-graphs can exist, which we do by construction. How would one block an edge? There are in principle two ways, adding and removing. With the additive approach one creates a two-step detour which forms a cycle once the edge in question is flipped. In the subtractive approach one merges unblocked edges with already blocked ones, decreasing the total number of (unblocked) edges. While the additive approach only threatens the finiteness of the constructed graph, the merging approach is more fickle: In order to glue two vertices together, they shall be distanced at least 4 edges apart, as otherwise one risks creating 1-, 2- or 3-cycles in the process. However, the constant presence of 3-cliques required for blocking edges densifies the graph enough that creating these distances is not trivial. The approach presented here consists of recursively and intricately merging many multi-pyramids together.

**Related Work.** There are several works on flip graphs of graph orientations. Most closely related, in [18] the authors show that for every $k \geq 2$ and for every $2k$-edge connected undirected graph $G$, any orientation of $G$ can be flipped to a $k$-edge connected orientation of $G$ via single orientation flips, where throughout the whole flip sequence the (oriented) connectivity is monotonically increasing. Moreover, they show how to find such a flip sequence in polynomial time. In a different direction, the work [3] presented at WG 2019 [2] studies flip graphs of so-called $\alpha$-orientations of a simple graph, in which every vertex $v$ has a specified outdegree $\alpha(v)$, and a flip consists of reversing all edges of a (directed) cycle. It is shown that deciding whether the flip distance is at most $k$ is NP-hard. They also use this result to show hardness of determining whether the flip distance between two perfect matchings of a graph $G$ is at most $k$ (where a flip exchanges edges and non-edges of the matching along an alternating cycle), even when $k = 2$ and $G$ is 2-connected, bipartite, subcubic, and planar.

In a more geometric context, flip graphs of different classes of plane straight-line graphs on a given finite point set in $\mathbb{R}^2$ have been extensively studied. This includes plane perfect matchings with the above-mentioned flip operation [13, 16, 24], as well as triangulations [22, 23, 19], spanning trees [6, 1, 7, 9, 14], spanning paths [4, 5, 10, 20, 26], perfect matchings [13, 16, 24], polygonizations [12], or partitions [17]. Further instances of flip graphs include combinatorial games, maybe most prominently the Rubik's cube [21, 27], but also flip graphs of non-geometric objects, such as bitstrings, permutations, combinations, and partitions [11]. For further examples of flip graphs and reconfiguration problems, as well as applications of them, we refer to the survey articles [15, 25, 8].

## 2  Formal Statement of the Main Conjecture

The following four definitions present no novel concepts, but instead serve to fix notation.

**Definition 1 (Graphs).** *A finite graph $G = (V, E)$ is a* simple oriented *graph iff it has no 1-, or 2-cycles (no loops or double edges), i.e. $E \subseteq \{(i, j) \in V \times V \mid i \neq j\}$ and $(i, j) \in E \implies (j, i) \notin E$*

For the remainder of this paper *graph* will refer to nonempty, finite, simple, oriented graphs.

**Definition 2 (3-Cycles and 3-Cliques).** *Let $G = (V, E)$ be a oriented graph and $a, b, c \in V$. The vertices $a, b, c$ form a* 3-cycle *iff there are edges $(a, b), (b, c)$ and $(c, a) \in E$. They form a* simplex *iff there are edges $(a, b), (b, c), (a, c) \in E$. In both cases they form a* 3-clique.

**Definition 3 (Vertex Identification).** *Let $G = (V, E)$ be a graph and $T \subseteq V$ non-empty. We denote the graph after identifying all vertices of $T$ to one new vertex with $G/T$. The edges are identified accordingly. Multiple identical (i.e. with the same orientation) edges between two vertices are regarded as one. If $T = \{a, b\}$ then we also write $a \sim b$.*

**Definition 4 (Edge-Flip).** *Let $G = (V, E)$ be a simple oriented graph. An* edge-flip *of $(i, j) \in E$ transforms the graph by replacing the edge $(i, j)$ with $(j, i)$.*

Note that [28] uses the notation SEF instead of *edge-flip*.

We state a slightly weaker version of the conjecture as presented in [28]. We drop the bound on number of edge-flips and the requirement for the final graph to be acyclic, focusing on 3-cycles instead. It is clear that refuting the following conjecture also refutes the one in [28].

**Conjecture 1.** *Does any finite, simple, oriented graph admit a sequence of edge-flips which monotonically decreases the number of 3-cycles to zero?*

## 3   Reduction to FBD-graphs

We show that Conjecture 1 may be falsified in the general case with, what we believe to be an easier concept, FBD-graphs.

The key idea exploited are blocked edges, see Figure 1. Note that in general, a blocked edge is not always *sufficiently* blocked to prohibit its inversion: there could be more 3-cycles associated to an edge than blocks.

**Definition 5 (Blocked Edge).** *Let $G = (V, E)$ be a simple, oriented graph. An edge $(a, b) \in E$ is* blocked *iff there exists a detour over vertex $c \in V$ such that $(a, c), (c, b) \in E$. An edge is* unblocked *iff no such detour exists. A graph is* fully blocked *iff all of its edges are blocked.*

In all subsequent figure blocked edges are *blue* whereas unblocked eddges are *red* assuming it is relevant. A vertex $v \in V$ is referred to as *blocked* iff all edges $\{(a, b) : (a, b) \in E, a = v \text{ or } b = v\}$ involving $v$ are blocked.

**Definition 6.** *An* FBD-graph *is a fully blocked, directed, finite, non-empty, simple, oriented graph without (directed) 3-cycles.*

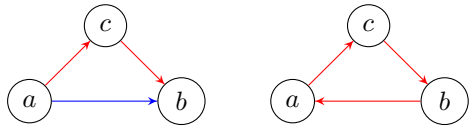

Fig. 1: Left: A blocked *blue* edge $(a,b)$ and its detour $(a,c),(c,b)$. Right: Flipping the blocked edge $(a,b)$ results in a 3-cycle. All unblocked edges are *red*

Verifying that a graph is FBD is easy, in particular with computer aid:

**Proposition 1.** *Let $G$ be a graph with adjacency matrix $A$. Then $G$ is an $FBD-graph$ if and only if:*

$$\begin{aligned}
\operatorname{tr} A &= 0, & simple \\
\operatorname{tr} A^2 &= 0, & oriented \\
\operatorname{tr} A^3 &= 0, & \text{no 3-cycles} \\
A^2 &\geq A, & \text{fully blocked}
\end{aligned}$$

*where the last $\geq$ is meant as an element-wise comparison.*

*Proof.* Recall that $A^i_{a,b}$ denotes the number of directed walks of length $i$ from $a$ to $b$. Thus $\operatorname{tr} A^i$ upper bounds the number of $i$-cycles, explaining the first three equations. Further, $A^2 \geq A$ states that for every edge $(a,b)$ there is at least one path of length 2 from $a$ to $b$ and hence every edge in the graph is blocked. □

The key reason to introduce FBD-graphs is that their existence contradicts Conjecture 1:

**Theorem 1.** *Let $D$ be an FBD-graph. Then there exists a simple, oriented graph $D'$ which does not admit to a sequence of edge-flips monotonically decreasing the number of 3-cycles to zero.*

*Proof.* Assume the existence of an FBD-graph $D = (V, E)$. We construct a new graph $D' = (V', E')$ using six copies of $D$, namely $D_i = (V_i, E_i), 1 \leq i \leq 6$. From every copy we pick one edge $(a_i, b_i) \in E_i$ and identify them the following way (see Figure 2):

$$d \triangleq a_1 \sim a_2 \sim b_5 \sim b_6, \quad e \triangleq a_3 \sim a_4 \sim b_1 \sim b_2, \quad f \triangleq a_5 \sim a_6 \sim b_3 \sim b_4.$$

This results in a 3-cycle $(d, e, f)$ of doubly blocked edges, i.e. edges that each create two new 3-cycles when flipped. Thus the 3-cycle $(d, e, f)$ is now sufficiently blocked: flipping any of its edges would increase the overall number of 3-cycles in $D'$. Therefore $D'$ provides a counterexample to Conjecture 1. □

The existence of any FBD-graph refutes Conjecture 1 in generality. This can not be easily extended to subclasses of graphs, except if the subclass is closed under above construction. Still, FBD-graphs appear to be the easier concept, as the following lemmata demonstrate:

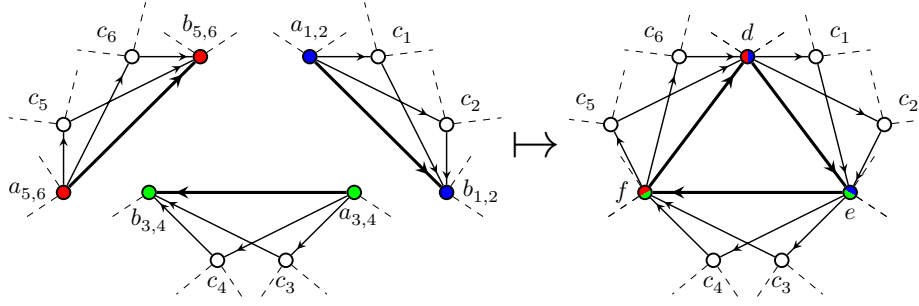

Fig. 2: Construction of the doubly blocked three-cycle in $D'$

**Lemma 1.** *There are no complete FBD-graphs.*

*Proof.* As the existence of a larger cycle in a complete graph implies the existence of 3-cycles, complete FBD-graph must be an acyclic complete graph. These express a total order on their vertices, and an edge between two neighbouring vertices w.r.t to this order cannot be blocked, as this would imply the existence of a vertex and thus value in between. Thus there are no complete FBD-graphs. □

This result does not follow immediately from the fact that Conjecture 1 holds true as shown in [28], as completeness is not closed under the construction of the reduction above. Still, it is noteworthy that this proof has a fraction of the complexity compared to the matching statement in [28].

The difference becomes even more pronounced with the following statement. Again, it does not immediately transfer to a statement that Conjecture 1 holds true for all graphs where every edge participates in at most 2 3-cycles. We believe that statement to be true, but at the time of writing a proof still eludes us.

**Lemma 2.** *Let $G = (V, E)$ and $e \in E$ and $\alpha(e)$ be defined as the number of 3-cycles $e$ participates in. If $\operatorname{avg}_{e \in E} \alpha < 3$, $G$ cannot be an FBD-graph. Similarly, if $\max_{e \in E} \alpha = 2$, $G$ cannot be an FBD-graph.*

*Proof.* A 3-clique consists of three edges, but can itself only block at most one edge. Thus for all edges to be blocked as required of an FBD-graph, each *blocked* edge must participate in at least three 3-cycles. But as in an FBD-graph all edges are blocked, this average extends to all edges, thus $\operatorname{avg}_{e \in E} \alpha \geq 3$. Of course $\max_{e \in E} \alpha$ must then also exceed 2. □

**Remark 1.** FBD-graphs appear to be the easier concept for two reasons: One no longer has to worry if a blocked edge is also *sufficiently* blocked, i.e. there are more blocks than 3-cycles connected to these edge, and they eliminate the potential for 3-cycles aimlessly wandering around, neither reducing the overall count of 3-cycles nor progressing in any way towards the goal.

While the first problem might cause only minor headaches, the second one is worth thinking about from a computational complexity angle: If $G$ is an FBD-graph or not is easily decidable in polynomial time (see Proposition 1). In contrast, deciding that $G$ is a counterexample to Conjecture 1 seems, on first glance, potentially NP-hard.

## 4 Constructing an FBD-graph

### 4.1 Overview

As shown in Theorem 1, constructing an FBD-graph is enough to refute Conjecture 1. We thus seek to construct one.

**The Construction in a Nutshell** Informally the procedure to construct an FBD-graph could be described as follows: First we construct a big patch of "inner" edges, which are all iteratively blocked except a "boundary" of unblocked edges. Then, by growing the patch large enough, we might be able to gain enough distance between the "innermost core" and the boundary that we may identify the (unblocked) boundary edges with the (blocked) inner edges, thus blocking the boundary in the process. This approach has a reasonable chance for success if we manage to achieve a minimal distance of 4 between the boundary vertices and the innermost vertices — any smaller distance might introduce 1-, 2- or 3-cycles during the identification, 4-cycles are fine however. One has to be careful to ensure that even after identifying a pair of vertices, the distance for other boundary vertices to the inner core has not shrunk below 4.

**The Construction in a Coconut-Shell** Instead of mere vertices and edges, our basic building block from now on is a multi-pyramid with a (oriented) 4-cycle as a base.

By nature of the multi-pyramid, all edges excluding the base are already blocked (See Figure 3). To block the base edges though, we glue (two) more multi-pyramids onto a previous one in a particular way (this procedure is described in detail in 4.2). Nothing is preventing us from blocking the base of those two newly introduced multi-pyramids, iterating upon this construction for as long as desired and allowing us to increase the distance from the center (original pyramid) to the boundary (the ones most recently introduced) to arbitrary numbers. For our specific counterexample we desire to iterate this procedure ten times (nine times is not enough).

This whole construction could be interpreted as a complete binary tree, each node representing a 4-cycle respectively a multi-pyramid, and two sibling nodes blocking all unblocked edges of their parents. From that perspective all layers except the last become the inner patch of blocked edges, and the last layer, the leaves, forms the unblocked boundary. It is also clear that the number of leaves and thus unblocked pyramid bases grows exponentially, in our case of ten levels leading to 1024 unblocked pyramids outnumbering the 1023 blocked ones on level 0-9. So far, this does not seem like a good premise to identify outer boundary pyramids with an inner core.

Luckily, from level 3 onwards, the leaves grow distant enough to (most) other leaves that we can identify them with each other, as long as we don't identify siblings (leaves sharing the same father) and first-degree cousins (leaves sharing the same grandfather). But cousins of 2nd or any higher degree are distant enough to be safely identified with each other without risking 1-, 2- or 3-cycles in the process, see 4.3. This observation reduces the number of unblocked pyramids from 1024 to a mere 4.

The last step is the long-foreshadowed fold of the boundary inwards: As the four remaining leaves on level 10 are distant enough to the four fully blocked nodes on level 2 (the distance grows by one for every two layers, see 4.4 for more details), we may now merge the these together, ridding the construction of all unblocked edges in the process. Using the four pyramids on level 2 as targets is minimal, as otherwise inter-leaf-distance would not be upheld in the identification process.

To formally verify that this process does not create any 1-, 2- or 3-cycles and that indeed every edge is now blocked we describe the graph in a computer-friendly way and check all FBD-graph-properties computationally, see 4.5.

### 4.2 Step I: Iteratively Growing a Patch of Blocked Pyramids

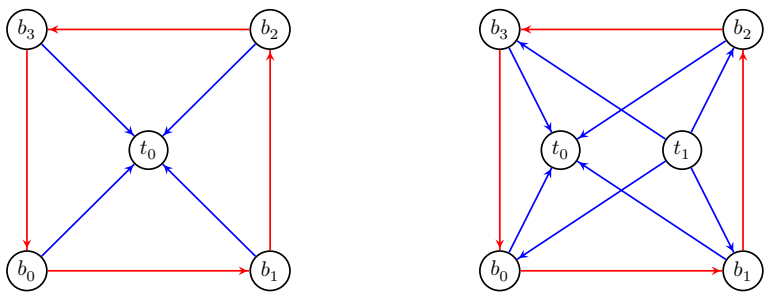

Fig. 3: A multi-pyramid with one tip on the left, and a tip and a pit to the right. Note that all edges towards/from tips and pits are blocked.

**Multi-Pyramids.** The *multi-pyramids* are the most essential building block to our construction. In fact, the finished FBD-graph consists of 2046 multi-pyramids intricately merged into each other. Our particular multi-pyramids consist of a *(pyramid) base* which is a 4-cycle, composed from four *base-edges*. With the base set, one may add arbitrarily many *tips*: every tip is connected via 4 *tip-edges* oriented towards the tip (forming a target) or towards the base (forming a source); see Figure 3. A neat property of pyramids is that *every* tip-edge is blocked, as the base-edges allow for detours. The base-edges remain unblocked.

To illustrate the method of iteratively blocking pyramids with other pyramids to increase the distance between the inner core and unblocked boundary we go through the process of creating a vertex ($t_c$) with distance of two (shorthanded $\text{dist}(t_c) = 2$) step-by-step. Everything that follows next is also depicted in Figure 4.

We start by defining the multi-pyramid $C = (V_C, E_C)$ with vertices $V_C = \{c_0, c_1, c_2, c_3, t_c\}$ with base-edges $E_C = \{(c_i, c_{(i+1) \bmod 4}) : i \in 0\ldots 3\}$ together with tip-edges $\{(c_i, t_c) : i \in 0\ldots 3\}$. The edges determining $\text{dist}(t_c) = 1$ are the edges of the base of $C$. To block the first two edges, i.e. $(c_1, c_2)$ and $(c_2, c_3)$, we first add an additional multi-pyramid $A$, where $A = (V_A, E_A)$ with $V_A = \{a_0, a_1, a_2, a_3, t_{at}, t_{as}\}$ and
$E_A = \{(a_i, a_{(i+1) \bmod 4}) : i \in 0\ldots 3\} \cup \{(t_{at}, a_i) : i \in 0\ldots 3\} \cup \{(a_i, t_{as}) : i \in 0\ldots 3\}$.
The two tips of the multi-pyramid $A$ are then merged with the vertices $c_1 \sim t_{at}, c_2 \sim a_0$ and $c_3 \sim t_{as}$. By this identification the two base-edges of $C$ have been upgraded to be tip-edges of the multi-pyramid $A$ and are thus blocked.

We repeat the same procedure for the other two base-edges $(c_0, c_1)$ and $(c_3, c_0)$ of $C$ by introducing another multi-pyramid $B$ (see Figure 4).

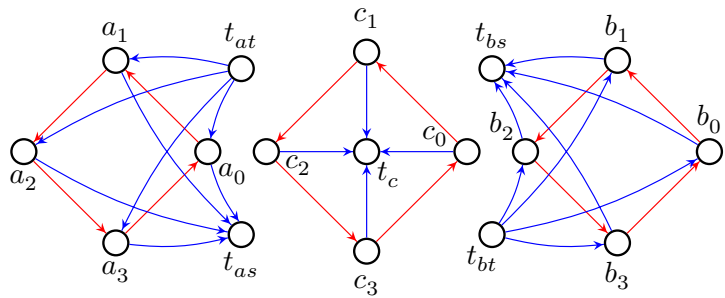

(a) Initially the 4-cycle $(c_0, c_1, c_2, c_3)$ is unblocked.

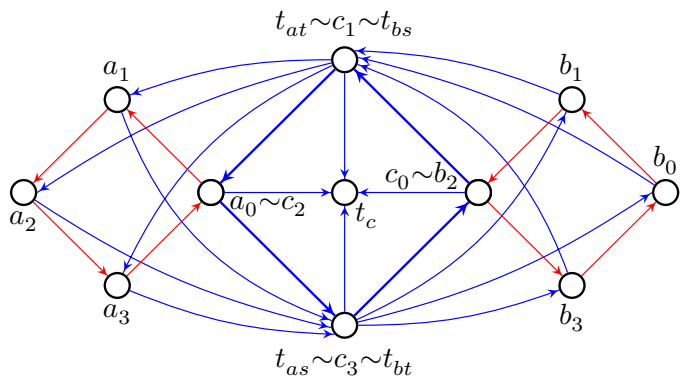

(b) After identification the 4-cycle $(c_0, c_1, c_2, c_3)$ is blocked.

Fig. 4: Blocking the base-edges of the initial pyramid $C$.

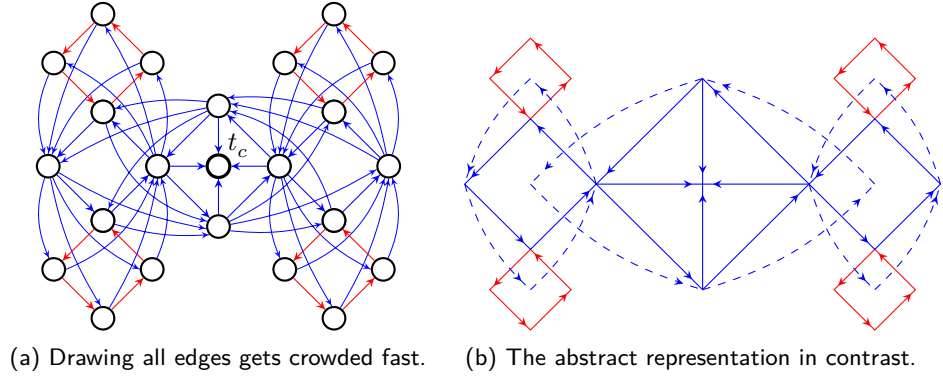

(a) Drawing all edges gets crowded fast.   (b) The abstract representation in contrast.

Fig. 5: Blocking the unblocked cycles a second time allows $\mathrm{dist}(t_c) = 2$.

While this procedure was successful in blocking all the base-edges of $C$ the still unblocked base-edges of $A$ (respectively $B$) cling to $c_2$ (respectively $c_2$), thus $\mathrm{dist}(t_c)$ remains 1. However, recursive repetition of the procedure mitigates this unfortunate effect. Figure 5 depicts the next iteration (and an abstraction thereof, which we will use in the future) in which the bases of $A$ and $B$ are blocked similar to the base of $C$, resulting in $\mathrm{dist}(t_c) = 2$.

There is no reason not to continue this pattern to reach any needed distance, except *potentially* the exponentially growing unblocked 4-cycles in the bases of the last iteration, i.e. $2^i$ after $i$ iterations.

### 4.3 Step II: Reducing the Boundary to Constant Size

Achieving arbitrary unblocked vertex distances at the expense of exponential number of unblocked edges seems like replacing one evil with another. We observe however that with more iterations of this procedure leaves grow more distant to most other leaves, and isn't a distance of 4 is enough to safely identify them with each other? Indeed, only one more iteration of aboves procedure is necessary to experience the beneficial effects of this phenomenon for the first time.

To observe this behaviour first hand and verify that all distance requirements are properly upheld, we iterate aboves construction three times, ending up with the graph depicted in Figure 6 on the left respectively a binary-tree abstraction on the right. We refer to the plain 4-cycle ($t_c$ is dropped from now on as it just served ilustratory purpose) in the center as to level 0. The bases of the multi-pyramids $A$ and $B$ lie on level 1, the bases of their four blocking pyramids on level 2 and so on. The eight pyramids of level 3 are split into two groups A and B, depending on which of the pyramids $A$ or $B$ is closer (see Figure 6).

We now verify that the minimal distance between any pairing of bases from A and B is four (see Figure 6), implying that we can boldly glue them together. More precisely, we glue all vertices of all four matched base-pairs together (the precise matching does not matter), without any risk of 1-, 2-, or 3-cycles. Note that this procedure works just as well with more iterations of the construction — the only leaves not mergeable are siblings and first-level-cousins. Thus, after the merger, only four leaves remain.

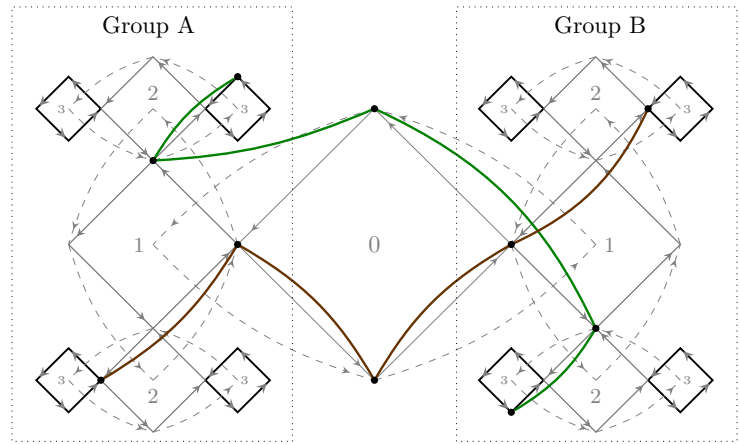

Fig. 6: The distance between vertices of level 3 4-cycles in Group A and Group B is at least 4. The numbers in the 4-cycles denote the respective level.

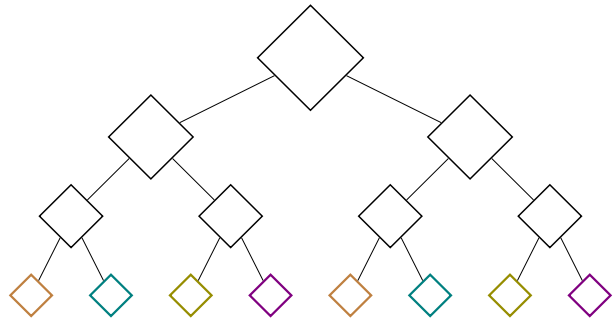

Fig. 7: The same graph as in Figure 6 interpreted as a binary tree. Nodes of the same color could now be identified with each other.

### 4.4 Step III: Folding the Remaining Boundary Edges Inwards

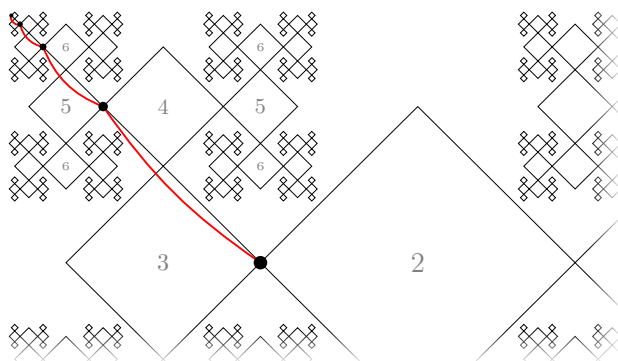

Fig. 8: 4-cycles of level 10 are of distance 4 to 4-cycles of level 2.

To get rid of the remaining 4 unblocked pyramids, we require a new tactic, foreshadowed in the overview: we shall merge the unblocked boundary with the innermost core.

Recalling that a minimal distance of 4 is safe for gluing parts of the graph together, and realizing that every two iterations of the construction represent a distance gain of one, on might guess that 8 levels are plenty to warrant enough distance between the innermost core (the center pyramid) and the leaves. However, we would violate post-merge inter-leaf-distance by merging them all onto one center-pyramid.

A simple solution is to pick four disjunct targets instead, which we get by extending the construction to level 10. Now we have four fully blocked pyramids on level 2 which are simultaneously far enough away from the leaves at level 2+8=10. Thus we can now merge the boundary with level 2. A very abstract representation is depicted in Figure 9.

An alternative but very similar construction would utilize 4 separate trees each of level 8 — this results in the same graph as removing level 0 and 1 (which now only serve illustratory purposes) from the finished construction.

### 4.5 Step IV: Computationally Verify FBD-properties

So far the entire construction avoided one crucial question: How do we ensure that merging pyramids of a certain distance does not critically decrease distances in future mergers?

As this question becomes rather messy, we decided to rely on a computer assisted proof, utilizing tools we already developed over the course of our investigation.

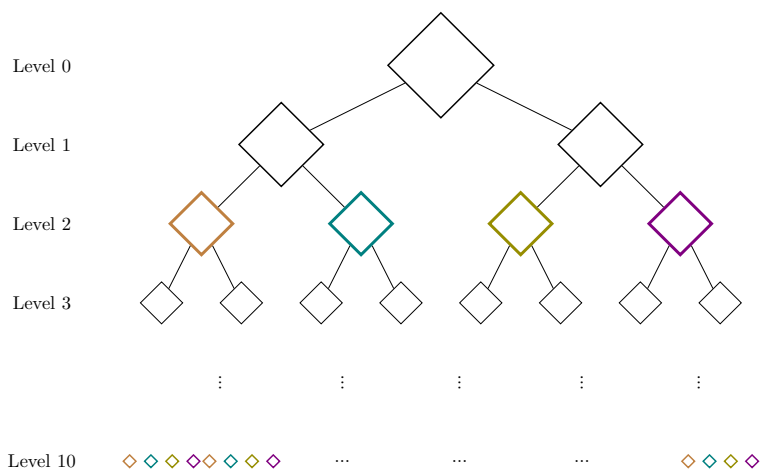

Fig. 9: The full construction in an abstract binary tree representation: Leaves on level 10 are identified with each other iff they match colors. These 4 disjunct groups are furthermore identified with the 4 distinct nodes on level 2.

In principle it would be enough to give one adjacency matrix representing an FBD-graph and then verify it via the (simple, linear algebra based) verification code that checks for all four conditions stated in Proposition 1. We encourage interested readers to do so themselves, the code and graphs (different varieties) are available at https://repository.tugraz.at/records/xfpsr-pzx32. Once run, it outputs:

```
$ curl -s https://repository.tugraz.at/records/xfpsr-pzx32/files/verifier.py \
> | python3
loaded 16368 edges
tr A¹ = 0
tr A² = 0
tr A³ = 0
A <= A²
all checks passed! graph is FBD
```

With this computational proof complete we now proudly claim:

### 4.6 Main Result

**Theorem 2.** *FBD-graphs exist.*

**Corollary 1.** *Conjecture 1 does not hold.*

## 5 Discussion and Future Work

We constructed a counterexample to Conjecture 1 and thus gave proof that the conjecture is not true in general. However, the construction presented here (12276 vertices and 98208 edges) is very large and very intricate, and does not leave much leeway for error. A reduction to slightly less than half the size is discussed in Section B.

Naturally one wonders if there are smaller FBD-graphs, as it helps in finding boundaries to the conjecture. Eventually, an entirely different approach to constructing FBD-graphs (if it exists) might be helpful.

We remark that the average number of 3-cliques per edge in the FBD-graph constructed in this work is three. In contrast, we show that graphs in which the maximum number of 3-cliques per edge is less than three cannot be FBD, indicating that we are at the theoretical limit.

This is however not true for the *maximal* number of 3-cliques per edge over the graph — here our construction requires 516, which we managed to reduce to 34 (see Section B). Can this be reduced to the theoretical lower limit of three? This immediately raises another elementary yet hard question: Do (undirected, nonempty) graphs where every edge participates in exactly three 3-cliques exist?[1]

In ongoing work, we investigate the relation between the conjecture and the existence of FBD-graphs from a complexity theoretic perspective, as already teased upon in Remark 1: What is the computational complexity of deciding whether a given graph forms a counterexample to Conjecture 1? We believe it to be significantly harder (NP-hard?) than deciding whether a graph is an FBD-graph, which is easily verifiable in polynomial time.

Finally, despite our work focusing on refuting Conjecture 1, there is reason to believe it actually holds true for a variety of graph classes: In [28] they already showed it for complete graphs, and the nonexistence of FBD-graphs where the maximal edge-3-clique participation rate is smaller than three indicates that the conjecture should be true for any graph with a maximal edge-3-clique participation smaller or equal to three. Actually, in the future full version of this paper, we provide a proof that Conjecture 1 holds true for any graph that admits an embedding on an oriented surface in which every 3-clique is the boundary of a face of the embedded graph. It would be interesting to see for which other graph classes the conjecture can be verified.

**Acknowledgements.** We would like to thank Michael Kerber for helpful discussions on the topic, and additionally Bettina Klinz, Yannic Maus and Jonathan Krebs for initial investigations on the problem.

---

[1] Chances are we are simply not aware of an already known solution.

## A Code to Verify FBD-graphs

### A.1 Concise Formalization of the FBD-graph

We've seen already that the iterative blocked pyramid construction can be interpreted as a tree: Take pyramid bases as nodes the nodes and the joint vertices (like $c_2$ and $a_0$ in Figure 4) as edges. This perspective allows for easy bookkeeping of all pyramid-bases and thus all the vertices, edges and identifications of our FBD-graph. Not only allows it for a concise description, it is greatly aids in writing graph-generating code, allowing the properties to be easily checked.

**Naming Scheme.** The root of this tree is the initial 4-cycle which we name with empty string $\varepsilon$. The bases of the multi-pyramids blocking $\varepsilon$ (i.e. $A$ and $B$) are named $\triangleleft$ and $\triangleright$. They themselves are blocked by $\triangleleft\triangleleft$, $\triangleleft\triangleright$ for $\triangleleft$ and analogously $\triangleright\triangleleft$ and $\triangleright\triangleright$ for $\triangleright$, and so on. Individual vertices are specified by appending $\alpha, \beta, \gamma$ or $\delta$ to the name of associated 4-cycle. Thus $\alpha$ is a vertex in level 0 and $\triangleleft\triangleleft\triangleright\gamma$ is a vertex in level 3. To refer to *all* bases of a certain level we write e.g. $\lozenge^{10}$ for all bases of level 10.

The bases of consecutive levels always share a vertex. While there is some leeway in specifying which vertex is merged with which, we use (in our code and the supplied resulting FBD-graph) the following identifications: A 4-cycle with name $x$ is blocked by identifying $x\beta$ with $x\triangleleft\alpha$ and $x\delta$ with $x\triangleright\alpha$ as illustrated in Figure 10. Further $x\alpha$ and $x\gamma$ are identified with the tips of $x\triangleleft$ and $x\triangleright$.

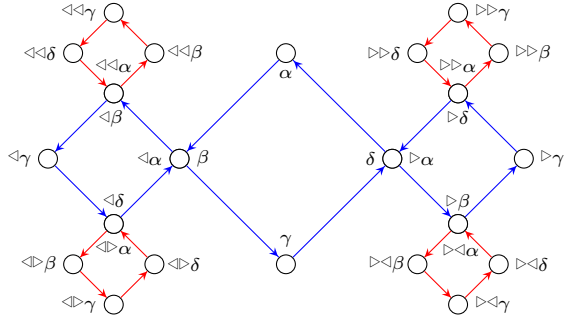

Fig. 10: Vertex naming scheme for the counterexample. For brevity only the base-edges of 3 levels are shown.

**All Folds in One Equation.** We build a recursively blocked pyramids ten levels deep. Then, with above naming scheme, we can now concisely describe all identifications via $\forall\ \Delta \in \{\triangleleft\triangleleft, \triangleleft\triangleright, \triangleright\triangleleft, \triangleright\triangleright\} \forall\ \Xi \in \{\alpha, \beta, \gamma, \delta\} : \lozenge^8 \Delta\Xi \sim \Delta\Xi$. The resulting graph is then carefully transformed into an adjacency matrix[2].

---

[2] The code for that is less curated and available on request. It might serve as a starting point to construct other FBD-graphs.

The following program then checks all conditions necessary for an adjancency matrix to represent an FBD-graph, it is available at https://repository.tugraz.at/records/xfpsr-pzx32.

```python
import numpy as np
from os import path
from urllib import request

fn = 'FBD_original.txt'         # original version as described in the paper
#fn = 'FBD_minimal_vertices.txt' # smaller graph w.r.t number of vertices/edges
#fn = 'FBD_minimal_alpha.txt'    # smaller graph w.r.t 3-clique factor alpha

# download graph data if not present
if not path.isfile(fn):
    request.urlretrieve(
        f'https://repository.tugraz.at/records/vg8bv-c8g87/files/{fn}',
        fn)

# load edges
edges = np.loadtxt(fn, dtype=int, delimiter=',')

# build adjacency matrix from edge list
a = np.zeros((edges.max()+1, edges.max()+1))
a[edges[:,0], edges[:,1]] = 1

print(f'loaded {np.sum(a, dtype=int)} edges')

# A^2, A^3
aa = a@a
aaa = aa@a

# tests
if np.trace(a) == 0:
    print(f'tr A¹ = 0')
else: raise Exception('graph is not simple')

if np.trace(aa) == 0:
    print(f'tr A² = 0')
else: raise Exception('graph is not oriented')

if np.trace(aaa) == 0:
    print(f'tr A³ = 0')
else: raise Exception('graph has 3-cycles')

if np.all(aa>=a):
    print(f'A <= A²')
else: raise Exception('graph contains unblocked edges')

print('all checks passed! graph is FBD')
```

Note that the smaller respectively less complex graphs briefly touched upon in Section 5 can be loaded and checked by modifying lines 5-7.

## B  Minimality of the Counterexample

Traditionally speaking, smaller graphs mean less vertices and edges. The construction presented in Section 4 has 2046 vertices and 16368 edges and is far from minimal, as it was instead optimized for simplicity. It is straightforward to see that 6 vertices and 48 edges can be trimmed off by removing the zeroth and first level of pyramids. But there is more potential to reduce: While pyramids are at most three steps apart, individual vertices are up to five steps apart. Identifying all vertices with distance 5 leaves us with 1888 vertices and 16256 edges, but remarkably some are now doubly blocked. This reduces the requirement of six copies of an FBD-graph as described in Theorem 1 to three, halving the size of the complete counterexample to Conjecture 1. Note that there are still still *more* vertices with distance 4 left and thus more potential for optimisation. Initial investigation convinced us that some systematic study would be beneficial.

Besides vertices and edges, we found another measurement of complexity intriguing: With $\alpha(e)$ given by the number of 3-cliques the edge $e$ participates in, we somewhat describe the local complexity. Any FBD-graph must have an $\operatorname{avg} \alpha \geq 3$: A 3-clique can block at most one edge, but is composed out of three edges, thus every blocked edge must on average participate in at least three 3-cliques. Our *original* construction hits $\operatorname{avg} \alpha = 3$ precisely, whereas the trimmed one has a slightly increased average as there are doubly blocked edges present. However, things become more interesting with $\max \alpha$, which with the construction in Section 4 takes a value of 516, as the base-edges of the pyramids on level 10 are highly interconnected. Although a slightly modified construction achieves to reduce $\max \alpha$ to 34, there is still a sizeable gap between the *current* theoretical lower limit of $\max \alpha = 3$.